\documentclass[11pt]{amsart}
\usepackage{amssymb}
\usepackage{amsmath}
\usepackage[dvips]{graphicx}
\usepackage{psfrag}
\begin{document}
\setlength{\parindent}{1.2em}
\def\COMMENT#1{}
\def\TASK#1{}
\def\noproof{{\unskip\nobreak\hfill\penalty50\hskip2em\hbox{}\nobreak\hfill%
       $\square$\parfillskip=0pt\finalhyphendemerits=0\par}\goodbreak}
\def\endproof{\noproof\bigskip}
\newdimen\margin   
\def\textno#1&#2\par{%
   \margin=\hsize
   \advance\margin by -4\parindent
          \setbox1=\hbox{\sl#1}%
   \ifdim\wd1 < \margin
      $$\box1\eqno#2$$%
   \else
      \bigbreak
      \hbox to \hsize{\indent$\vcenter{\advance\hsize by -3\parindent
      \sl\noindent#1}\hfil#2$}%
      \bigbreak
   \fi}
\def\proof{\removelastskip\penalty55\medskip\noindent{\bf Proof. }}
\def\enddiscard{}
\long\def\discard#1\enddiscard{}
\newtheorem{firstthm}{Proposition}
\newtheorem{thm}[firstthm]{Theorem}
\newtheorem{prop}[firstthm]{Proposition}
\newtheorem{lemma}[firstthm]{Lemma}
\newtheorem{cor}[firstthm]{Corollary}
\newtheorem{problem}[firstthm]{Problem}
\newtheorem{defin}[firstthm]{Definition}
\newtheorem{conj}[firstthm]{Conjecture}
\newtheorem{theorem}[firstthm]{Theorem}
\newtheorem{claim}{Claim}
\def\eps{{\varepsilon}}
\def\N{\mathbb{N}}
\def\R{\mathbb{R}}
\def\P{\mathcal{P}}
\def\Ybar{\bar{Y}}

\title{$k$-ordered Hamilton cycles in digraphs}
\author{Daniela K\"uhn,  Deryk Osthus and Andrew Young }
\date{}
\maketitle \vspace{-.8cm}
\begin{abstract} \noindent
Given a digraph~$D$, let $\delta^0(D):=\min\{\delta^+(D), \delta^-(D)\}$ be the
minimum semi-degree of~$D$. $D$ is $k$-ordered Hamiltonian if 
for every sequence $s_1,\dots,s_k$ of distinct vertices of~$D$ there is a directed Hamilton cycle
which encounters $s_1,\dots,s_k$ in this order. 
Our main result is that every digraph~$D$ of sufficiently large order~$n$
with $\delta^0(D)\ge \lceil(n+k)/2\rceil -1$ is $k$-ordered Hamiltonian. The bound on the
minimum semi-degree is best possible. An undirected version of this result was proved earlier by 
Kierstead, S\'ark\"ozy and Selkow~\cite{KSSordered}.
\end{abstract}
\section{Introduction}

The famous theorem of Dirac determines the smallest minimum degree of a graph
which guarantees the existence of a Hamilton cycle.
There are many subsequent results which investigate degree conditions that 
guarantee the existence of a Hamilton cycle with some additional properties.
In particular, Chartrand (see~\cite{NgSchultz}) introduced the notion of a 
Hamilton cycle which has to visit a given set of vertices in a prescribed order.
More formally, we say that a graph $G$ is \emph{$k$-ordered} if 
for every sequence $s_1,\dots,s_k$ of distinct vertices of~$G$ there is a cycle
which encounters $s_1,\dots,s_k$ in this order. $G$~is \emph{$k$-ordered Hamiltonian}
if it contains a Hamilton cycle with this property. 
Kierstead, S\'ark\"ozy and Selkow~\cite{KSSordered} showed that for all $k \ge 2$ 
every graph on $n \ge 11k-3$ vertices of minimum degree at least $\lceil n/2 \rceil + \lfloor k/2 \rfloor-1$
is $k$-ordered Hamiltonian. This bound on the minimum degree is best possible and proved
a conjecture of Ng and Schultz~\cite{NgSchultz}.
Several related problems have subsequently been considered:
for instance, the case when $k$ is large compared to $n$ was investigated 
in~\cite{Faudreeetal} (but has not been completely settled yet).
Ore-type conditions were investigated in~\cite{NgSchultz,Faudreeetal,CGPf}.
For more results in this direction, see the survey by Gould~\cite{Gould}.

It seems that digraphs provide an equally natural setting for such problems.
Our main result is a version of the result in~\cite{KSSordered} for digraphs.
The digraphs we consider do not have loops and we allow at most one edge in each
direction between any pair of vertices.
Given a digraph~$D$, the \emph{minimum semi-degree~$\delta^0(D)$} of~$D$ is the
minimum of the minimum outdegree~$\delta^+(D)$ of~$D$ and its minimum
indegree~$\delta^-(D)$.
\begin{thm}\label{thm:ordHam}
For every $k \ge 3$ there is an integer $n_0=n_0(k)$ such that every 
digraph~$D$ on~$n\ge n_0$ vertices with $\delta^0(D)\ge \lceil(n+k)/2\rceil -1$
is $k$-ordered Hamiltonian.
\end{thm}
Our proof shows that one can take $n_0:=Ck^9$ where $C$ is a sufficiently large constant.
Note that if~$n$ is even and~$k$ is odd the bound on the minimum semi-degree is
slightly larger than in the undirected case.
However, it is best possible in all cases. In fact, if the minimum semi-degree is smaller, 
it turns out that $D$ need not even be $k$-ordered. This is easy to see if $k$ is even: 
let~$D$ be the digraph which consists of a complete digraph~$A$ of order~$\lceil n/2\rceil +k/2-1$
and a complete digraph~$B$ of order~$\lfloor n/2\rfloor +k/2$ which has precisely $k-1$ vertices in
common with~$A$. Pick vertices $s_1,s_3,\dots,s_{k-1}\in A-B$ and $s_2,s_4,\dots,s_{k}\in B-A$.
Then~$D$ has no cycle which encounters $s_1,\dots,s_k$ in this order.
A similar construction also works if both $k$ and $n$ are odd.
The construction in the remaining case is a little more involved, see~\cite{KOdilinked} for details.
Note that every Hamiltonian digraph is $2$-ordered Hamiltonian, so the case when $k \le 2$ in Theorem~\ref{thm:ordHam}
is covered by the result of Ghouila-Houri~\cite{GhouilaHouri} (Theorem~\ref{thm:GH} below)
which implies that every digraph with
minimum semi-degree at least $n/2$ contains a Hamilton cycle.

Theorem~\ref{thm:ordHam} can be used to deduce a version for edges which have to be traversed in a prescribed 
order by the Hamilton cycle: we say that a digraph~$D$ is \emph{$k$-arc ordered Hamiltonian} if, for
every sequence $e_1,\dots,e_k$ of independent edges, $D$ contains a Hamilton cycle
which encounters $e_1,\dots,e_k$ in this order. $D$ is \emph{$k$-arc Hamiltonian} if
it contains a Hamilton cycle which encounters these edges in any order.
$D$ is called \emph{Hamiltonian $k$-linked} if $|D|\ge 2k$
and if for every sequence $x_1,\dots,x_k,y_1,\dots,y_k$ of distinct vertices there are disjoint paths
$P_1,\dots,P_k$ in~$D$ such that $P_i$ joins~$x_i$ to~$y_i$ and such that together all the~$P_i$ cover
all the vertices of~$D$. Thus every digraph~$D$ which is Hamiltonian $k$-linked is also
$k$-arc ordered Hamiltonian. Indeed, if $x_1y_1,\dots,x_ky_k$ are the (directed) edges our
Hamilton cycle has to encounter then disjoint paths linking~$y_{i-1}$ to~$x_i$ for all $i=1,\dots,k$
yield the required Hamilton cycle.
\begin{cor}\label{thm:arcordHam}
For all $k \ge 3$ there is an integer $n_0=n_0(k)$ such that every 
digraph~$D$ on~$n\ge n_0$ vertices with $\delta^0(D)\ge \lceil n/2\rceil +k-1$
is Hamiltonian $k$-linked and thus in particular $k$-arc ordered Hamiltonian.
\end{cor}
The examples in~\cite{KOdilinked} show that in both parts of Corollary~\ref{thm:arcordHam}
the bound on the minimum semi-degree is best possible. In fact, if the minimum semi-degree
is smaller then one cannot even guarantee the digraph to be $k$-arc ordered.%
      \COMMENT{For Hamiltonian $k$-linkedness the bound is best possible as it is
best possible for just $k$-linkedness. This automatically implies that it is best possible for
``$k$-arc ordered Hamiltonian''. One way to see this is to note that
the minimum degree which guarantees a digraph to be $k$-arc ordered Hamiltonian
already guarantees it to be $k$-linked. Indeed, if we wish to link the pairs $x_i,y_i$
then add the edges $y_{i-1}x_i$ to $D$ to obtain a digraph $D'$ (if $y_{i-1}x_i\notin D$).
Then $\delta^0(D')\ge \delta^0(D)\ge \lceil n/2\rceil +k-1$ and so
$D'$ is $k$-arc ordered Hamiltonian. In particular, $D'$ has a Hamilton cycle
containing the edges $y_kx_1,y_1x_2,\dots,y_{k-1}x_k$ in this order. The paths
between these edges yield the required spanning set of $x_i$-$y_i$ paths in~$D$.
Another way to see that the bound on the min degree is best possible for 
``$k$-arc ordered Hamiltonian'' is to look at the following extremal examples.
If $n$ is even let~$D$ consist of complete digraphs~$A$ and~$B$
of order~$n/2 +k-1$ which have precisely $2k-2$ vertices $x_1,\dots,x_{k-1},y_1,\dots,y_{k-2},y_k$
in common. Let $x_k$ be a vertex in $A\setminus B$ and let $y_{k-1}$ be a vertex in $B\setminus A$.
Then $x_1y_1,\dots,x_ky_k$ are edges but no cycle contains them in this order.
So suppose next that~$n$ is odd. In this case, we define~$D$ as follows.
Let~$A$ and~$B$ be disjoint complete digraphs of order $\lceil n/2\rceil-k-1$.
Add a complete digraph~$X$ of order~$2k-3$ with vertices $x,x_1,x_4,\dots,x_k,y_3,\dots,y_k$
and join all vertices in~$X$ to all vertices in~$A\cup B$ with edges in both directions.
Add a set $S:=\{x_2,x_3,y_1,y_2\}$ of~4 new vertices such that each vertex in~$S$
is joined to each vertex in~$X$ with edges in both directions.
Moreover, we add all the edges between different vertices in~$S$ except
for~$y_1x_2$ and~$y_2x_3$. Finally, we connect the vertices in~$S$ to the vertices
in~$A\cup B$ as follows. Both~$y_1$ and~$x_2$ receive edges from every vertex in~$B$ and send
edges to every vertex in~$A$. Additionally, $y_1$ will receive an edge from every vertex in~$A$
and $x_2$ will send an edge to every vertex in~$B$.
Both~$y_2$ and~$x_3$ receive edges from every vertex in~$A$ and send
edges to every vertex in~$B$. Additionally, $y_2$ will receive an edge from every vertex in~$B$
and $x_3$ will send an edge to every vertex in~$A$. This is the same example as in the linkedness-paper
but our $y_1$ is the $x_1$ there, our $x_2$ is the $y_1$ there, our $y_2$ is the $x_2$ there
and our $x_3$ is the $y_2$ there. So in particular, the minimum degree is the right one.
To see that~$D$ is not $k$-arc-ordered note that $x_1y_1,\dots,x_ky_k$ are edges
but no cycle contains them in this order as this would mean that we could link
$y_i$ to $x_{i+1}$ by disjoint paths.} 
A result of Bermond~\cite{Bermond} (see also~\cite{Berge})
implies that if $\delta^0(D)\ge \lceil (n+k)/2\rceil$
then~$D$ is $k$-arc Hamiltonian. It easily follows that if
$\delta^0(D)\ge \lceil (n+1)/2\rceil$, then $D$ is Hamiltonian $1$-linked,
i.e.~Hamiltonian connected (see~\cite{Berge}).%
     \COMMENT{If $x$ and $y$ are the vertices we wish to link then add the edge $yx$ to
$D$ to obtain a digraph $D'$ (if $yx\notin D$). Then $\delta^0(D')\ge \lceil (n+1)/2\rceil$ and so
$D'$ is $1$-arc Hamiltonian. In particular, $D'$ has a Hamilton cycle
containing the edge $yx$, this yields a Hamilton $x$-$y$ path in~$D$.}
This covers the case $k=1$ of Corollary~\ref{thm:arcordHam}.
As observed in~\cite[Thm~9.2.10]{digraphsbook}, if $\delta^0(D)\ge \lceil n/2\rceil+1$, then $D$ is 
Hamiltonian $2$-linked, which covers the case $k=2$ of Corollary~\ref{thm:arcordHam}. 
     
Corollary~\ref{thm:arcordHam} can easily be deduced from Theorem~\ref{thm:ordHam}
as follows:
let $x_1,\dots,x_k$ and $y_1,\dots,y_k$ be distinct vertices where we aim to link
$x_i$ to $y_i$ for all $i$. 
Let~$D'$ be the digraph obtained from~$D$ by contracting~$x_i$
and~$y_{i-1}$ into a new vertex~$s_i$ whose outneighbourhood is that of~$x_i$
and whose inneighbourhood is that of~$y_{i-1}$. More precisely,
let $A:=\{x_1,\dots,x_k,y_1,\dots,y_k\}$. Then~$D'$ is
the digraph obtained from $D-A$ by adding new vertices $s_1,\dots,s_k$
and defining the edges incident to these new vertices as follows.
The outneighbours of~$s_i$ are the outneighbours of~$x_i$ in $V(D)\setminus A$ as well as
all the~$s_j$ for all those $j\neq i-1$ for which~$y_j$ is an outneighbour of~$x_i$ in~$D$ (where
$y_0:=y_k$). Similarly, inneighbours of~$s_i$ are the inneighbours of~$y_{i-1}$ in $V(D)\setminus A$
as well as all the~$s_j$ for all those $j\neq i$ for which~$x_j$
is an inneighbour of~$y_{i-1}$ in~$D$.  It is easy to check that%
      \COMMENT{Indeed, if $x\in D-A$ then e.g. $d^+_{D'}(x)\ge d^+_D(x)-k
\ge \lceil n/2\rceil +k-1-k=\lceil (|D'|+k)/2\rceil -1$. The calculation
for~$s_i$ is similar. For example, $d^+_{D'}(s_i)\ge d^+_D(x_i)-k$ as we loose
(at most) all the $x_j$ with $j\neq i$ as well as $y_{i-1}$.}
$\delta^0(D')\ge \lceil (|D'|+k)/2\rceil -1$ and that a Hamilton cycle in~$D'$
which encounters $s_1,\dots,s_k$ in this order corresponds to 
a spanning set of disjoint paths from~$x_i$ to~$y_i$. 

A result of Chen et al.~\cite[Theorem~10]{Chenetalk-arc} implies
that the smallest minimum degree which guarantees an undirected graph to be $k$-arc
ordered Hamiltonian is $\lfloor n/2\rfloor +k-1$.%
      \COMMENT{The graph consisting of two cliques whose orders are as equal as possible
and which have $2k-2$ vertices in common shows that their result is best possible.
The (oriented) edges are chosen as in the directed case.}
(A graph is $k$-arc ordered Hamiltonian
if for any sequence of~$k$ independent oriented edges there exists a Hamilton cycle which encounters
these edges in the given order and orientation.)
The smallest minimum degree which forces a graph to be $k$-linked was determined
by Kawarabayashi, Kostochka and Yu~\cite{KKY}. It is not clear whether the minimum
degree for Hamiltonian $k$-linkedness is the same.%
      \COMMENT{One cannot apply the trick above to the result in~\cite{KSSordered}.}

The main tool in our proof of Theorem~\ref{thm:ordHam} is a recent result by the 
first authors (Theorem~\ref{thm:ord} below), 
which shows that the degree condition in Theorem~\ref{thm:ordHam} at least
guarantees a $k$-ordered cycle
(but not necessarily a Hamiltonian one). The strategy of the proof of Theorem~\ref{thm:ordHam}
is to consider such a cycle of maximal length and to show that it must be Hamiltonian. 
The same strategy was already applied in the proof of the undirected case in~\cite{KSSordered}.
However, both parts of the strategy are more difficult in the digraph case: the existence of a
$k$-ordered directed cycle (i.e.~Theorem~\ref{thm:ord}) already confirms a conjecture
of Manoussakis~\cite{Manoussakis} for large~$n$.
The Hamiltonicity of a $k$-ordered cycle of maximal length is easier to show in the undirected
case as one can consider `local transformations' of a given $k$-ordered cycle which reverse the orientation of certain 
segments of the cycle. This means that apart from some basic observations like Lemma~\ref{longloop}
below our proof is quite different from that in~\cite{KSSordered}. 
\begin{thm}{\cite{KOdilinked}} \label{thm:ord}
Let~$k$ and~$n$ be integers such that $k\ge 2$ and~$n\ge 200k^3$.
Then every digraph~$D$ on~$n$ vertices with $\delta^0(D)\ge \lceil(n+k)/2\rceil -1$
is $k$-ordered.
\end{thm}

\section{Notation and tools}
Given a digraph~$D$, we write~$V(D)$ for its vertex set, $E(D)$ for its edge set
and~$|D|:=|V(D)|$ for its order. 
We write $xy$ for the edge directed from~$x$ to~$y$. More generally, if~$A$ and~$B$
are disjoint sets of vertices of~$D$ then an~\emph{$A$-$B$ edge} is an edge of the
form~$ab$ where $a\in A$ and $b\in B$. A digraph is \emph{complete} if every pair of distinct
vertices is joined by edges in both directions.

Given disjoint subdigraphs~$D_1$ and~$D_2$ of a digraph~$D$ such that $D_1\cup D_2$ is spanning
and a set $A\subseteq V(D_1)$, we write~$N^+_{D_i}(A)$ for the set of all those
vertices $x\in V(D_i)\setminus A$ which in the digraph~$D$ receive an edge from some vertex in~$A$.
$N^-_{D_i}(A)$ is defined similarly. If~$A$ consists of a single vertex~$x$,
we just write $N^+_{D_i}(x)$ etc.~and put $d^+_{D_i}(x):=|N^+_{D_i}(x)|$ and
$d^-_{D_i}(x):=|N^-_{D_i}(x)|$. So in particular, $N^+_D(x)$ is the outneighbourhood
of~$x$ in~$D$ and $d^+_{D}(x)$ is its outdegree. Also, note that $N^+_{D_1}(x)$ is the
outneighbourhood of~$x$ in the subdigraph~$D[V(D_1)]$ of~$D$ induced by~$V(D_1)$
and not its outneighbourhood in~$D_1$
(where $x\in D_1$). We let $N_D(x):=N^+_D(x)\cup N^-_D(x)$.

If we refer to paths and cycles in digraphs then we always mean that
they are directed without mentioning this explicitly. The \emph{length} of a path is the
number of its edges. Given two vertices $x,y\in D$,
an  \emph{$x$-$y$ path} is a path which is directed from~$x$ to~$y$.
Given two vertices~$x$ and~$y$ on a directed cycle~$C$, we write $xCy$ for the subpath of~$C$
from $x$ to~$y$. Similarly, given two vertices~$x$ and~$y$ on a directed path~$P$
such that~$x$ precedes~$y$, we write $xPy$ for the subpath of~$P$ from $x$ to~$y$.

A digraph~$D$ is \emph{strongly connected} if for every ordered pair~$x,y$ of vertices of~$D$ there
exists an $x$-$y$ path. $D$ is \emph{Hamiltonian connected} if for every ordered pair~$x,y$
of vertices of~$D$ there exists a Hamilton path from~$x$ to~$y$. (So Hamiltonian connectedness
is the same as Hamiltonian 1-linkedness.)

We will often use the following result of Ghouila-Houri~\cite{GhouilaHouri} which gives a sufficient condition
for the existence of a Hamilton cycle in a digraph. In particular, it implies
a version of Theorem~\ref{thm:ordHam} for $k\le 2$ as any Hamiltonian digraph is 2-ordered Hamiltonian.

\begin{thm}\label{thm:GH}
Suppose that~$D$ is a strongly connected digraph such that $d^+_D(x)+d^-_D(x)\ge |D|$
for every vertex $x\in D$. Then~$D$ is Hamiltonian.
\end{thm}

The next result of Overbeck-Larisch~\cite{OL} provides
a sufficient condition for a digraph to be Hamiltonian connected.

\begin{thm}\label{thm:Hamconn}
Suppose that~$D$ is a digraph such that $d^+_D(x)+d^-_D(y)\ge |D|+1$ whenever~$xy$
is not an edge. Then~$D$ is Hamiltonian connected.
\end{thm}

\section{Preliminary results}\label{sec:prelim}

Let~$D$ be a digraph satisfying the conditions of Theorem~\ref{thm:ordHam}.
Let $S=(s_1,\dots,s_k)$ by any sequence of~$k\ge 3$ vertices of~$D$.
We will often view~$S$ as a set. An \emph{$S$-cycle} in~$D$ is a cycle which encounters
$s_1,\dots,s_k$ in this order. So we have to show that~$D$ has a Hamiltonian $S$-cycle.
Theorem~\ref{thm:ord} implies the existence of an~$S$-cycle in~$D$.
Let~$C$ be a longest such cycle and suppose that~$C$ is not Hamiltonian.
Let~$H$ be the subdigraph of~$D$ induced by all the vertices outside~$C$.
Our aim is to find a longer $S$-cycle by modifying~$C$ (yielding a contradiction). 
The purpose of this section is to collect the properties of~$C$ and~$H$ that
we need in our proof of Theorem~\ref{thm:ordHam}. 

We let~$F$ be the set of all those vertices on~$C$ which receive an edge
\emph{from} some vertex in~$H$ and we let~$T$ be the set of all those vertices
on~$C$ which send an edge \emph{to} some vertex in~$H$. Given~$i\in\N$,
we write~$F_i$ for the set of all those vertices on~$C$ which receive an edge
from at least~$i$ vertices in~$H$. Thus $F_1=F$. $T_i$ is defined similarly.
Given a vertex~$x$ on~$C$,
we will denote its successor on~$C$ by~$x^+$ and its predecessor by~$x^-$.

\begin{lemma}\label{Hlemma1}
$H$ is Hamiltonian connected and $d^-_H(x)+d^+_H(y)\ge |H|+k-2$
for all vertices $x,y\in H$. Moreover any digraph obtained from~$H$
by deleting at most~$2$ vertices is strongly connected and
$k\le |H|\le \lfloor \frac{n-k}{2}\rfloor$.
\end{lemma}
\proof
We first show that any two (not necessarily distinct) vertices~$x,y\in H$ for which~$H$ contains an
$x$-$y$ path, $P$ say, satisfy the degree condition in the lemma.
To see this, note that no vertex in~$N^-_C(x)$ is a predecessor of
some vertex in~$N^+_C(y)$. Indeed, if $v\in N^-_C(x)$ and $v^+\in N^+_C(y)$
then by replacing the edge~$vv^+$ with the path $vxPyv^+$
we obtain a longer $S$-cycle, a contradiction. But this means that
$d^-_C(x)+d^+_C(y)\le |C|$ and thus
\begin{equation}\label{eq:H2}
d^-_H(x)+d^+_H(y)\ge 2\left(\left\lceil \frac{n+k}{2}\right\rceil-1\right)-|C|
\ge |H|+k-2,
\end{equation}
as required. However, as $k\ge 3$ this degree condition means that
$N^-_H(x)\cap N^+_H(y)\neq \emptyset$ and so~$H$ contains an
$y$-$x$ path of length~2. Thus whenever~$H$ contains an $x$-$y$ path it also
contains a $y$-$x$ path.

Now let~$x$ and~$z$ be any two vertices of~$H$. What we have shown above
applied with $y:=x$ implies that $d^-_H(x)+d^+_H(x)\ge |H|+1$ and thus $|N_H(x)|\ge (|H|+1)/2$.
Note that by the above~$x$ is joined to every vertex in~$N_H(x)$ with paths in both directions. 
Similarly, $|N_H(z)|\ge (|H|+1)/2$ and~$z$ is joined to every vertex in~$N_H(z)$ with
paths in both directions. As $|N_H(x)\cap N_H(z)|>0$ this means that~$x$
is joined to~$z$ with paths in both directions, i.e.~$H$ is strongly connected.
Together with~(\ref{eq:H2}) this in turn implies that
$d^-_H(x)+d^+_H(z)\ge |H|+k-2\ge |H|+1$ for all vertices $x,z\in H$.
In particular,~$H$ is Hamiltonian connected by Theorem~\ref{thm:Hamconn}.

To show that any digraph $H'$ obtained from~$H$ by deleting at most~2 vertices
is strongly connected note that $d^-_{H'}(x)+d^+_{H'}(y)\ge |H'|-1$
for every~$x,y\in H'$. Thus if $x\neq y$ then either~$yx$ is an edge or~$H'$ contains an
$y$-$x$ path of length~2.

It now remains to prove the bounds on~$|H|$. Consider any vertex $x\in H$.
Then $2(|H|-1)\ge d^-_H(x)+d^+_H(x)\ge |H|+k-2$ and so $|H|\ge k$.
For the upper bound, note that no vertex in~$T$ has a successor in~$F$.
Indeed, if~$v$ is such a vertex in~$T$ and~$v^+$ is its successor then we could
replace~$vv^+$ with a path through~$H$ to obtain a longer $S$-cycle,
a contradiction. But this means that some vertex of~$C$ must have all its
inneighbours on~$C$ or all its outneighbours on~$C$. Thus $|C|\ge \lceil (n+k)/2\rceil$
and so $|H|\le \lfloor (n-k)/2\rfloor$.
\endproof

Recall that the proof of Lemma~\ref{Hlemma1} implies the following.

\begin{cor}\label{FTcor}
No vertex on~$C$ which lies in~$T$ has a successor in~$F$.
\end{cor}

The next result deals with the case when the vertices $x_1\in T$ and $x_2\in F$ are
further apart.

\begin{lemma}\label{longloop}
Suppose that~$x_1, x_2\in C$ are distinct and the interior of $x_1Cx_2$ does not contain a
vertex from~$S$. Then there are no distinct vertices $y_1,y_2\in H$ such that
$x_1y_1, y_2x_2\in E(D)$.
\end{lemma}
\proof
Suppose that such~$y_1,y_2$ do exist. Furthermore, we may assume that~$x_1$ and~$x_2$
are chosen such that they
satisfy all these properties and subject to this $|x_1Cx_2|$ is minimum.
Let~$Q$ denote the set of all vertices in the interior
of $x_1Cx_2$. Then our choice of~$x_1$ and~$x_2$ implies that
$N^-_C(y_1)\cap Q=\emptyset$ and $N^+_C(y_2)\cap Q=\emptyset$. 
Moreover, by Corollary~\ref{FTcor} no vertex in~$N^-_C(y_1)$ is a predecessor of
some vertex in~$N^+_C(y_2)$. Thus $d^-_C(y_1)+d^+_C(y_2)\le |C|-|Q|+1$
and so
$$
n+k-2\le d^-_D(y_1)+d^+_D(y_2) \le |C|-|Q|+1 +2(|H|-1)=n-|Q|+|H|-1.
$$
This implies that $|H|>|Q|$ and thus replacing the interior of $x_1Cx_2$ with a Hamilton path from~$y_1$
to~$y_2$ through~$H$ (which exists by Lemma~\ref{Hlemma1})
yields a longer $S$-cycle, a contradiction.
\endproof

The next two results will be used in the proof of Lemma~\ref{H2paths}.

\begin{lemma}\label{strongconn}
Let~$G$ be a digraph such that $d^+_G(x)+d^-_G(x)\ge |G|+3$ for every vertex $x\in G$
and $d^+_G(x)+d^-_G(y)\ge |G|+1$ for every pair of vertices $x,y\in G$.
Let~$z_1$ and~$z_2$ be distinct vertices of~$G$ such that $z_1z_2\notin E(G)$.
Then there exists a vertex $a\in N^+_G(z_1)\cap N^-_G(z_2)$ such that
$G-\{z_1,z_2,a\}$ is strongly connected.
\end{lemma}
\proof
First note that $|N^+_G(z_1)\cap N^-_G(z_2)|\ge 3$ since $z_1z_2\notin E(G)$. Pick
$a_1,a_2,a_3\in N^+_G(z_1)\cap N^-_G(z_2)$. We will show that one of these~$a_i$
can play the role of~$a$. Let~$G^*:=G-\{z_1,z_2\}$. Note that
$d^+_{G^*}(x)+d^-_{G^*}(x)\ge |G^*|+1$ for every vertex $x\in G^*$
and $d^+_{G^*}(x)+d^-_{G^*}(y)\ge |G^*|-1$ for every pair of vertices $x,y\in G^*$.
In particular, the latter condition implies that~$G^*$ is strongly connected.
Thus~$G^*$ has a Hamilton cycle~$C$ by Theorem~\ref{thm:GH}. Let~$a^+_1$ denote the
successor of~$a_1$ on~$C$ and let~$a^-_1$ be its predecessor.
Put $N^+:=N^+_{G^*}(a^-_1)\setminus \{a_1\}$ and $N^-:=N^-_{G^*}(a^+_1)\setminus \{a_1\}$.
Note that $|N^+|, |N^-|\ge 1$ since $d^+_{G^*}(a^-_1)+d^-_{G^*}(a^-_1)\ge |G^*|+1$
and $d^+_{G^*}(a^+_1)+d^-_{G^*}(a^+_1)\ge |G^*|+1$. Similarly $|N^+|+|N^-|\ge |G^*|-3$.
Clearly, if $a^-_1a^+_1$ is an edge or $N^+\cap N^-\neq \emptyset$, then~$G^*-a_1$ is strongly
connected and so we can take~$a$ to be~$a_1$. So we may assume that neither of these is the
case. But then $N^+\cup N^-=V(G^*)\setminus\{a_1,a^+_1, a^-_1\}$. Let $v\in N^+$ be such that
$|vCa^-_1|$ is maximal. Similarly, let $w\in N^-$ be such that $|a^+_1Cw|$ is maximal.
Note that if $w\in vCa^-_1$ then~$G^*-a_1$ is strongly connected.
So we may assume that this is not the case. But then~$v$ must be the successor of~$w$ on~$C$,
$N^+$ must consist of precisely the vertices in $V(vCa^-_1)\setminus \{a^-_1\}$ and $N^-$ must consist of
precisely the vertices in $V(a^+_1Cw)\setminus \{a^+_1\}$.

Let~$A^+:=N^+\cup \{a^-_1\}$ and $A^-:=N^-\cup \{a^+_1\}$. We may assume that~$G$ does
not contain an $A^+$-$A^-$ edge as otherwise~$G^*-a_1$ is strongly connected.
We will now show~$G^*[A^+]$ is complete and that~$a_1$ receives an edge from every
vertex in~$A^+$.
So consider any vertex $x\in A^+$. Then $d^+_{G^*}(x)+d^-_{G^*}(a^+_1)\ge |G^*|-1$.
Together with the fact that there is no $A^+$-$A^-$ edge this shows that
$N^+_{G^*}(x)=(A^+\cup \{a_1\})\setminus \{x\}$. Thus~$G^*[A^+]$ is complete
and~$a_1$ receives an edge from every vertex in~$A^+$.
Similarly one can show that~$G^*[A^-]$ is complete and that~$a_1$ sends an edge to every
vertex in~$A^-$.

Now consider~$a_2$ and~$a_3$. If for example $a_2\neq v,w$ then~$G^*-a_2$ is strongly connected
and so we can take~$a$ to be~$a_2$.  As one can argue similarly for~$a_3$, we may assume
that $v=a_2$ and $w=a_3$.
If $a^+_1a^-_1$ is an edge or $a_1\in N^+_{G^*}(a^+_1)\cap N^-_{G^*}(a^-_1)$
then~$G^*-a_2$ is strongly connected. (Here we used that $a^-_1\neq v=a_2$ since
$|N^+|\ge 1$.)
If this is not the case, then $d^+_{G^*}(a^+_1)+ d^-_{G^*}(a^-_1)\ge |G^*|-1$ implies
the existence of some vertex $x\in N^+_{G^*}(a^+_1)\cap N^-_{G^*}(a^-_1)$ with $x\neq a_1$.
If $x\in A^+$ then $a^+_1x$
is an $A^-$-$A^+$ edge avoiding~$w=a_3$ and so~$G^*-a_3$ is strongly connected.
(Here we used that $a^+_1\neq w=a_3$ since $|N^-|\ge 1$.) Similarly, if
$x\in A^-$ then~$G^*-a_2$ is strongly connected. 
Altogether, this shows that we can take~$a$ to be~$a_1$, $a_2$ or~$a_3$.
\endproof 

\begin{lemma}\label{lowdegH}
Suppose that~$H$ contains a vertex~$v$ with $d^-_H(v)+d^+_H(v)\le |H|+k-1$.
Suppose that $x_1,x_2\in T$ and~$y_1,y_2\in F$ are distinct vertices on~$C$.
Then $x_1v, vy_1\in E(D)$ or $x_2v, vy_2\in E(D)$
(or both).
\end{lemma}
\proof
Let~$F_v$ denote the set of all those vertices on~$C$ which receive an edge from~$v$.
Let~$T^+_v$ denote the set of all those vertices on~$C$ whose predecessor sends an
edge to~$v$. Corollary~\ref{FTcor} implies that $T^+_v\cap F_v=\emptyset$. 
Since
$$
d^-_C(v)+d^+_C(v)\ge 2\left(\left\lceil \frac{n+k}{2}\right\rceil-1\right)-(|H|+k-1)
\ge |C|-1
$$
this shows that at most one vertex on~$C$ lies outside $T^+_v\cup F_v$.
Let~$z$ be the vertex in $V(C)\setminus (T^+_v\cup F_v)$ (if it exists).

Suppose first that $z\notin F$ (this also covers the case when~$z$ does not exist).
Then $z\neq y_1,y_2$. Also either
$z\neq x^+_1$ or $z\neq x^+_2$. So let us assume that $z\neq x^+_1$ (the case when $z\neq x^+_2$
is similar). We will show that
$x_1v, vy_1\in E(D)$. So suppose first that $x_1v\notin E(D)$.
Then $x^+_1\notin T^+_v$ and thus $x^+_1\in F_v$, a contradiction to Corollary~\ref{FTcor}.
Similarly, if $vy_1\notin E(D)$ then $y_1\notin F_v$ and thus $y_1\in T^+_v$,
i.e.~the predecessor of~$y_1$ lies in~$T$, contradicting Corollary~\ref{FTcor}.

So suppose next that $z\in F$ and thus, by Corollary~\ref{FTcor}, the predecessor of~$z$
does not lie in~$T$. This in turn implies that $z\neq x^+_1,x^+_2$.
Moreover either $z\neq y_1$ or $z\neq y_2$. So let us assume that $z\neq y_1$.
Similarly as before one can show that $x_1v,vy_1\in E(D)$.%
       \COMMENT{Suppose first that $x_1v\notin E(D)$. Then 
$x^+_1\notin T^+_v$ and thus $x^+_1\in F_v$ (as $z\neq x^+_1$), a contradiction to
Corollary~\ref{FTcor}. Similarly, if $vy_1\notin E(D)$ then
$y_1\notin F_v$ and thus $y_1\in T^+_v$, i.e.~the predecessor of~$y_1$ lies
in~$T$, contradicting Corollary~\ref{FTcor}.}
\endproof

In our proof of Theorem~\ref{thm:ordHam} we will frequently need two disjoint
paths through~$H$ joining two given disjoint pairs of vertices on~$C$ in order
to modify~$C$ into a longer $S$-cycle. The following lemma implies the existence
of such paths provided that the pairs consist of vertices having sufficiently many
neighbours in~$H$ (see also Corollary~\ref{corH2paths}).

\begin{lemma}\label{H2paths}
Suppose that $X_1,X_2\subseteq T$ and $Y_1,Y_2\subseteq F$ are
disjoint subsets of~$V(C)$ such that
$|N^+_H(X_1)|, |N^+_H(X_2)|\ge 3$ and $|N^-_H(Y_1)|, |N^-_H(Y_2)|\ge 3$.
Then there are disjoint $X_i$-$Y_i$ paths~$P_i$ of length at least~$2$
and such that all inner vertices of~$P_1$ and~$P_2$ lie in~$H$.
Moreover, if $|H|\ge 15$ and if we even have that
$|N^+_H(X_1)|, |N^+_H(X_2)|\ge 8$ and $|N^-_H(Y_1)|, |N^-_H(Y_2)|\ge 8$ then we can find such paths
which additionally satisfy $|P_1\cup P_2|\ge |H|/6$.
\end{lemma}
\proof By disregarding some neighbours if necessary we may assume that $|N^+_H(X_1)|=|N^+_H(X_2)|=|N^-_H(Y_1)|=|N^-_H(Y_2)|$.
Our first aim is to show that for some $i\in \{1,2\}$ there is an $X_i$-$Y_i$ path~$P_i$
which satisfies the following properties:
\begin{enumerate}
\item[(i)] The graph~$H':=H-V(P_i)$ has a Hamilton cycle~$C'$.
\item[(ii)] All $x,y\in H'$ satisfy $d^+_{H'}(x)+d^-_{H'}(y)\ge |H'|-2$.
\item[(iii)] $3\le |P_i|\le 5$, i.e.~$P_i$ contains at least~1 and at most~3 vertices
from~$H$.
\item[(iv)] If $i=1$ then $|N^+_H(X_2)\cap V(P_1)|\le 2$ and $|N^-_H(Y_2)\cap V(P_1)|\le 2$.
If $i=2$ then $|N^+_H(X_1)\cap V(P_2)|\le 2$ and $|N^-_H(Y_1)\cap V(P_2)|\le 2$.
\end{enumerate}
If we have found such an~$i$, say $i=1$, then our aim is to use the Hamilton cycle~$C'$
in order to find~$P_2$.
To prove the existence of such an~$i$, recall that Lemma~\ref{Hlemma1}
implies $d^-_H(x)+d^+_H(y)\ge |H|+k-2\ge |H|+1$ for every pair of vertices $x,y\in H$.
Thus condition~(ii) will hold automatically if~(iii) holds.

Now suppose first that there exists a vertex $z_1\in N^+_H(X_1)\cap N^-_H(Y_1)$.
Take~$i=1$ and take~$P_1$ to be any $X_1$-$Y_1$~path whose interior consists precisely
of~$z_1$. Then $d^-_{H'}(x)+d^+_{H'}(x)\ge |H'|$ for every $x\in H'$.
As~$H'$ is strongly connected by Lemma~\ref{Hlemma1} we can apply Theorem~\ref{thm:GH}
to find a Hamilton cycle~$C'$ of~$H'$. (If $|H'|=2$ then~$C'$ will consist of just a
double edge.) In the case when $N^+_H(X_2)\cap N^-_H(Y_2)\neq \emptyset$ we
proceed similarly.

Now suppose that $N^+_H(X_1)\cap N^-_H(Y_1)= \emptyset$ and
$N^+_H(X_2)\cap N^-_H(Y_2)= \emptyset$. Then Lemma~\ref{lowdegH}
implies that $d^-_H(x)+d^+_H(x)\ge |H|+k\ge |H|+3$ for every $x\in H$.
If there is an $N^+_H(X_1)$-$N^-_H(Y_1)$ edge~$z_1z_2$ take~$i:=1$ and take~$P_1$ to be
any $X_1$-$Y_1$ path whose interior consists of this edge. Then
$d^-_{H'}(x)+d^+_{H'}(x)\ge |H|-1= |H'|+1$ for every $x\in H'$ and so again, as~$H'$ is
strongly connected by Lemma~\ref{Hlemma1}, we can apply Theorem~\ref{thm:GH}
to find a Hamilton cycle~$C'$ of~$H'$. In the case when there is an $N^+_H(X_2)$-$N^-_H(Y_2)$
edge we proceed similarly.

Thus we may assume that $N^+_H(X_i)\cap N^-_H(Y_i)=\emptyset$ and that there is no
$N^+_H(X_i)$-$N^-_H(Y_i)$ edge (for $i=1,2$). Pick any vertex $z_1\in N^+_H(X_1)$
and let $z_2\in N^-_H(Y_1)$ be a vertex such that $|N^+_H(X_2)\cap \{z_1,z_2\}|\le 1$ and
$|N^-_H(Y_2)\cap \{z_1,z_2\}|\le 1$. (The fact that we can choose such a~$z_2$ follows
from $N^+_H(X_i)\cap N^-_H(Y_i)=\emptyset$ and our assumption that the sizes of the
$N^+_H(X_i)$ and the $N^-_H(Y_i)$ are equal.)%
       \COMMENT{Indeed, suppose for example that $z_1\in N^+_H(X_2)$ (and thus
$z_1\notin N^-_H(Y_2)$). So we can take $z_2$ to be any vertex inside
$N^-_H(Y_1)\setminus N^+_H(X_2)$. But $z_1\notin N^-_H(Y_1)$ since
$N^+_H(X_1)\cap N^-_H(Y_1)=\emptyset$ and so $z_1\in N^+_H(X_2)\setminus N^-_H(Y_1)$.
But this implies that $N^-_H(Y_1)\setminus N^+_H(X_2)$ is nonempty as
$|N^-_H(Y_1)|= |N^+_H(X_2)|$.}%
Apply Lemma~\ref{strongconn} with $G:=H$ to find a vertex
$z_3\in N^+_H(z_1)\cap N^-_H(z_2)$ such that $H-\{z_1,z_2,z_3\}$ is strongly connected.
Take $i:=1$ and~$P_1$ to be any $X_1$-$Y_1$ path whose interior consists of $z_1z_3z_2$.
Then $d^-_{H'}(x)+d^+_{H'}(x)\ge |H|-3= |H'|$ for every $x\in H'$
and so again~$H'$ contains a Hamilton cycle~$C'$ by Theorem~\ref{thm:GH}.
Our choice of~$z_1$ and~$z_2$ implies that~(iv) holds.

Altogether, this shows that in each case for some~$i$ there exists a path~$P_i$ satisfying (i)--(iv).
We may assume that~$i=1$. As mentioned before, our aim now is to use the Hamilton cycle~$C'$
of~$H'$ in order to find an $X_2$-$Y_2$ path~$P_2$ through~$H'$. In the case when
$|N^+_H(X_2)|,|N^-_H(Y_2)|\ge 3$ this is trivial
since by~(iv) both $N^+_H(X_2)$ and $N^-_H(Y_2)$ meet~$H'$ in at least one vertex.

So suppose now that $|H|\ge 15$ and $|N^+_H(X_2)|,|N^-_H(Y_2)|\ge 8$ and thus we wish to find
a long $X_2$-$Y_2$ path. To do this, let $N^+:=N^+_H(X_2)\cap V(H')$ and
$N^-:=  N^-_H(Y_2)\cap V(H')$. Thus $|N^+|,|N^-|\ge 6$ by~(iv).
Choose $a_1\in N^+$ and $b_1\in N^-$ to be distinct such that $|a_1C'b_1|$ is maximum.
If $|a_1C'b_1|\ge |H'|/6$ then we can take~$P_2$ to be any $X_2$-$Y_2$~path whose
interior consists of~$a_1C'b_1$. So we may assume that $|a_1C'b_1|\le |H'|/6$.

Note that the choice of~$a_1$ and~$b_1$ implies that $N^+,N^-\subseteq V(a_1C'b_1)$.
Moreover, all the vertices in~$N^+$ must precede the vertices in~$N^-$ on $a_1C'b_1$.
(Indeed, if e.g.~$a\in N^+$ and $b\in N^-$ are distinct vertices such that~$b$
precedes~$a$, i.e.~$a$ lies on~$bC'b_1$ then $|aC'b|\ge |H'|-|a_1C'b_1|\ge |H'|/2$,
contradicting the choice of~$a_1$ and~$b_1$.)
Thus $|N^+\cap N^-|\le 1$ and there are vertices $a_2,\dots,a_5\in N^+$ and
$b_2,\dots,b_5\in N^-$ such that $a_1,\dots,a_5,b_5,\dots,b_1$ are distinct and
appear on~$C'$ in this order. We now distinguish several cases.

\medskip

\noindent\textbf{Case~1.} \emph{There are $i,j\le 4$ such that $a_ib_j$ is an edge.}

\smallskip

\begin{figure}[htb!]
\begin{center}\footnotesize
\psfrag{1}[][]{$a_1$}
\psfrag{2}[][]{$a_2$}
\psfrag{3}[][]{$a_5$}
\psfrag{4}[][]{$b_5$}
\psfrag{5}[][]{$b_2$}
\psfrag{6}[][]{$b_1$}
\psfrag{7}[][]{$P_2$}
\psfrag{8}[][]{$v$}
\includegraphics[scale=0.5]{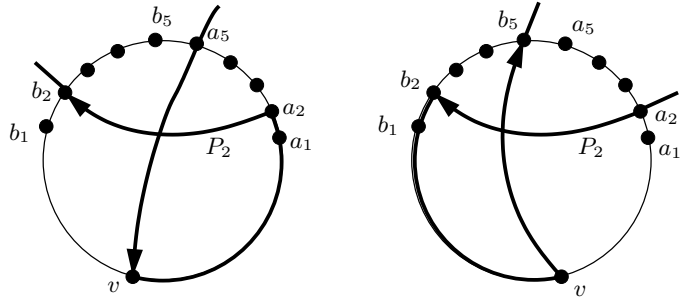}  
\caption{The path~$P_2$ in Case~1. The left figure is for the subcase when
$d^+_{H'}(a_5)\ge |H'|/2-1$ and the right figure is for the subcase when
$d^-_{H'}(b_5)\ge |H'|/2-1$.}
\label{fig:H2edge}
\end{center}
\end{figure}

\noindent
Note that $d^+_{H'}(a_5)\ge |H'|/2-1$ or $d^-_{H'}(b_5)\ge |H'|/2-1$ by~(ii). Suppose that
the former holds (the other case is similar). As $|a_1C'b_1|\le |H'|/6$
this means that~$a_5$ has at least $|H'|/3-1$ outneighbours in the interior of $b_1C'a_1$ and so
we can find such an outneighbour~$v$ with $|vC'a_1|\ge |H'|/3$.
But then we can take~$P_2$ to be any $X_2$-$Y_2$~path whose interior consists of
$a_5vC'a_ib_j$ (Figure~\ref{fig:H2edge}).

\medskip

\noindent\textbf{Case~2.} \emph{For all $i,j\le 4$ $a_ib_j$ is not an edge.}

\medskip

\noindent\textbf{Case~2.1.} \emph{There exists some vertex
$u\in N^+_{H'}(a_1)\cap N^-_{H'}(b_3)$.}

\smallskip

\noindent
Note that $u\neq a_2,b_4$ since by our assumption neither~$a_2b_3$ nor~$a_1b_4$ is an edge.
As before, either $d^+_{H'}(a_2)\ge |H'|/2-1$ or $d^-_{H'}(b_4)\ge |H'|/2-1$.
Suppose that the former holds (the other case is similar).

If~$u$ lies in the interior of~$a_1C'b_3$, let~$v$ be an outneighbour of~$a_2$
in the interior of~$b_3C'a_1$ with $|vC'a_1|\ge |H'|/3$. Then we can take~$P_2$ to be
any $X_2$-$Y_2$~path whose interior consists of $a_2vC'a_1ub_3$.

\begin{figure}[htb!]
\begin{center}\footnotesize
\psfrag{1}[][]{$a_1$}
\psfrag{2}[][]{$a_2$}
\psfrag{3}[][]{$b_3$}
\psfrag{4}[][]{$v$}
\psfrag{5}[][]{$u$}
\psfrag{6}[][]{$P_2$}
\psfrag{7}[][]{$b_1$}
\psfrag{8}[][]{$b_5$}
\includegraphics[scale=0.55]{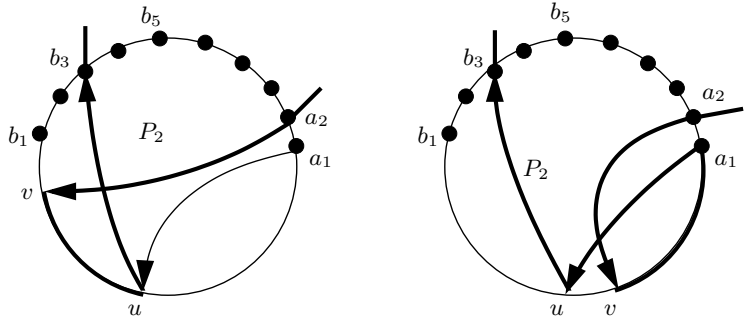}  
\caption{The path~$P_2$ in Case~2.1 if~$u$ lies in the interior of~$b_3C'a_1$.
The left figure is for the subcase when the interior of $b_3C'u$ contains at least $|H'|/6-1$
outneighbours of~$a_2$. The right figure is for the subcase when the interior
of $uC'a_1$ contains at least $|H'|/6-1$ outneighbours of~$a_2$.}
\label{fig:H2case22}
\end{center}
\end{figure}

So we may assume that~$u$ lies in the interior of~$b_3C'a_1$.
But then either the interior of $b_3C'u$ contains at least $|H'|/6-1$
outneighbours of~$a_2$ or the interior of $uC'a_1$ contains at least $|H'|/6-1$
outneighbours of~$a_2$. If the former holds let~$v$ be any outneighbour of~$a_2$
in the interior of $b_3C'u$ such that $|vC'u|\ge |H'|/6$ and take~$P_2$
to be any $X_2$-$Y_2$~path whose interior consists of $a_2vC'ub_3$
(see Figure~\ref{fig:H2case22}).
If the latter holds let~$v$ be any outneighbour of~$a_2$
in the interior of $uC'a_1$ such that $|vC'a_1|\ge |H'|/6$ and take~$P_2$
to be any $X_2$-$Y_2$~path whose interior consists of $a_2vC'a_1ub_3$.

\medskip

\noindent\textbf{Case~2.2.} \emph{There exists some vertex
$u\in N^+_{H'}(a_3)\cap N^-_{H'}(b_1)$.}

\smallskip

\noindent
This case is similar to Case~2.1 and we omit the details.

\medskip

\noindent\textbf{Case~2.3.} \emph{Both $N^+_{H'}(a_1)\cap N^-_{H'}(b_3)$
and $N^+_{H'}(a_3)\cap N^-_{H'}(b_1)$ are empty.}

\smallskip

\noindent
Together with~(ii) and our assumption that~$a_1b_3$ is not an edge this implies
that $N^+_{H'}(a_1)\cup N^-_{H'}(b_3)=V(H')\setminus \{a_1,b_3\}$.
Since~$a_3b_3$ is not an edge this means that~$a_1a_3$ is an edge.
Similarly it follows that~$b_3b_1$ is an edge. But as before either
$d^+_{H'}(a_2)\ge |H'|/2-1$ or $d^-_{H'}(b_2)\ge |H'|/2-1$. Suppose that
the former holds (the other case is similar). Then we can find an outneighbour~$v$
of~$a_2$ in the interior of $b_1C'a_1$ with $|vC'a_1|\ge |H'|/3$.
But then we can take~$P_2$ to be any $X_2$-$Y_2$~path whose interior consists of
$a_2vC'a_1a_3C'b_1$.
\endproof

Lemma~\ref{H2paths} immediately implies the following corollary, which is sometimes
more convenient to apply.

\begin{cor}\label{corH2paths}
Suppose that $x_1,x_2\subseteq T_3$ and $y_1,y_2\subseteq F_3$ are
distinct vertices on~$C$. Then~$D$ contains disjoint $x_i$-$y_i$ paths~$P_i$ of
length at least~$2$ such that all inner vertices of~$P_1$ and~$P_2$ lie in~$H$.
Moreover, if $|H|\ge 15$ and if we even have that $x_1,x_2\subseteq T_8$ and
$y_1,y_2\subseteq F_8$ then we can find such paths
which additionally satisfy $|P_1\cup P_2|\ge |H|/6$.
\end{cor}

The last of our preliminary results gives a lower bound
on the sizes of~$T_3$ and~$F_3$.

\begin{lemma}\label{sizeT5F5}
We have that $|T|,|F|\ge (n+k)/2-|H|$. Moreover, $|T_3|,|F_3|\ge (n-k)/2-|H|$ 
and $|T_3\cup F_3|\ge |C|-|H|-2k$.
\end{lemma}
\proof To see the bound on~$|T|$, note that $d^-_C(x)\ge \delta^0(D)-(|H|-1)
\ge (n+k)/2-|H|$ for every vertex $x\in H$ and so $|T|\ge (n+k)/2-|H|$.
The proof for~$|F|$ is similar.
To prove the bound on~$|T_3|$, we double-count the number $e(T,H)$ of
edges in~$D$ from~$T$ to~$V(H)$. Since $d^-_C(x)\ge(n+k)/2-|H|$ for
any vertex $x\in H$ we have that $e(T,H)\ge |H|((n+k)/2-|H|)$. 
On the other hand $e(T,H)\le |T_3||H|+2(|T|-|T_3|)=|T_3|(|H|-2)+2|T|$.
Before we can use this to estimate~$|T^3|$, we need an upper bound on~$|T|$.
For this, recall that $|F|\ge (n+k)/2-|H|$. Together with
Corollary~\ref{FTcor} this shows that $|T|\le |C|-|F|\le (n-k)/2$.
Altogether this gives
\begin{align*}
|T_3| & \ge \frac{|H|((n+k)/2-|H|)-(n-k)}{|H|-2}=
\frac{(|H|-2)(n-k)/2 -|H|(|H|-k)}{|H|-2}\\
& \ge \frac{n-k}{2}-|H|
= \frac{|C|-|H|-k}{2}.
\end{align*}
The proof for~$|F_3|$ is similar.
The bound on $|T_3\cup F_3|$ follows since $|T_3\cap F_3|\le k$.
Indeed, the latter holds since Lemma~\ref{longloop}
implies that whenever $s,s'\in S$ are distinct and no vertex from~$S$ lies in
the interior of $sCs'$ then $T_3\cap F_3$ meets $sCs'$ in at most one vertex.
\endproof

\section{Proof of Theorem~\ref{thm:ordHam}}

Throughout this section, we assume that the order~$n$ of our given digraph~$D$ is
sufficiently large compared to~$k$ for our estimates to hold. We will also omit floors and
ceilings whenever this does not affect the argument.
Let~$S$, $C$ and~$H$ be as defined at the beginning of Section~\ref{sec:prelim}.
Recall that we assume that~$C$ is not Hamiltonian and will show that we can
extend~$C$ into a longer $S$-ordered cycle (which would yield a contradiction
and thus would prove Theorem~\ref{thm:ordHam}).
Given consecutive vertices~$s,s'\in S$, we call the path obtained from~$sCs'$ by deleting~$s'$
the \emph{interval from~$s$ to~$s'$}. Thus no vertex from~$S$ lies in the interior of $sCs'$
and~$C$ consists of precisely~$|S|=k$ disjoint intervals.
In our proof of Theorem~\ref{thm:ordHam} we distinguish the following
4~cases according to the order of~$H$. Recall that~$|H|\ge k$ by Lemma~\ref{Hlemma1}.

\medskip

\noindent\textbf{Case~1.} $k\le |H|\le 220k^3$.

\smallskip

\noindent
Recall that $|T_3|\ge (n-k)/2-|H|\ge n/3$ by Lemma~\ref{sizeT5F5}
and so at least one of the~$k$ intervals of~$C$ must contain at least $n/(3k)$
vertices from~$T_3$. Suppose that this is the case for the interval~$I$ from~$s$
to~$s'$. Recall that by Lemma~\ref{sizeT5F5} at most~$|H|+2k\le 3|H|$ vertices of~$C$
do not lie in~$T_3\cup F_3$ and by Corollary~\ref{FTcor}
no vertex in~$F_3$ is the successor of a vertex in~$T_3$. Since every maximal subpath of~$I$
consisting of vertices from~$T_3$ is succeeded by at least one vertex outside~$T_3\cup F_3$,
it follows that~$I$ contains a subpath~$A$ which consists entirely of vertices from~$T_3$
and satisfies $|A|\ge n/(3k(3|H|+1))$. Let~$A_1$ be the subpath of~$A$ consisting of
its initial $n/(20k|H|)$ inner vertices and let~$A_2$ be the subpath of~$A$ consisting of
its last $n/(20k|H|)$ inner vertices.

Let~$t$ be the first vertex of~$A$. (So~$t^+$ is the first vertex of~$A_1$.)
Consider any vertex~$a$ on~$t^+Cs'$. Lemma~\ref{longloop} implies that $a\notin F$.
Thus $N^-_D(a)\subseteq V(C)$
and hence
\begin{equation}\label{eq:dega}
d^-_C(a)\ge \delta^0(D)\ge (n+k)/2-1\ge n+k-1-|H|-|F|>|C|-|F|.
\end{equation}
(To see the third inequality recall that $|F|\ge (n+k)/2-|H|$ by Lemma~\ref{sizeT5F5}.)

\medskip

\noindent\textbf{Case~1.1.} \emph{There are vertices $a_1\in A_1$ and $a_2\in A_2$
such that~$a_1a_2$ is an edge.} 

\smallskip

\begin{figure}[htb!]
\begin{center}\footnotesize
\psfrag{1}[][]{$a_1$}
\psfrag{2}[][]{$a^+_1$}
\psfrag{3}[][]{$a^-_2$}
\psfrag{4}[][]{$a_2$}
\psfrag{5}[][]{$w$}
\psfrag{6}[][]{$w^+$}
\psfrag{7}[][]{$P$}
\psfrag{8}[][]{$C$}
\psfrag{9}[][]{$a^-_1$}
\psfrag{10}[][]{$a_1$}
\psfrag{11}[][]{$a^-_2$}
\psfrag{12}[][]{$a_2$}
\psfrag{13}[][]{$w_1$}
\psfrag{14}[][]{$w^+_1$}
\psfrag{15}[][]{$w_2$}
\psfrag{16}[][]{$w^+_2$}
\psfrag{17}[][]{$P_1$}
\psfrag{18}[][]{$P_2$}
\includegraphics[scale=0.55]{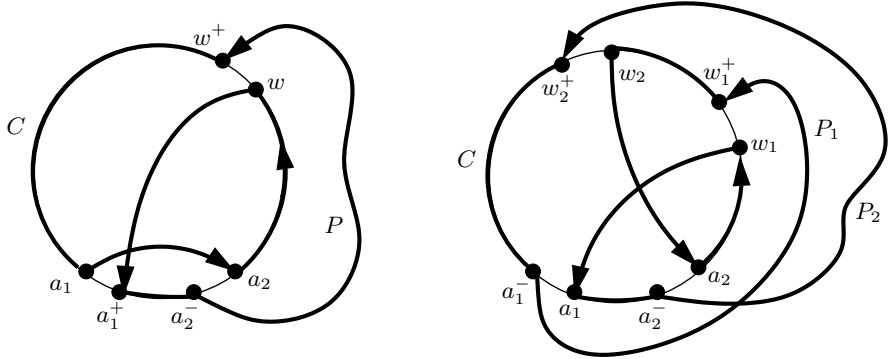}  
\caption{Extending~$C$ into a longer $S$-ordered cycle in Case~1.1 (left) and
Case~1.2 (right).}
\label{fig:proofcase1}
\end{center}
\end{figure}

\noindent
Inequality~(\ref{eq:dega}) applied with $a:=a^+_1$ 
implies that there exists a vertex $w\in N^-_C(a^+_1)$ such that the successor~$w^+$
of~$w$ lies in~$F$. Recall that~$F$ avoids~$t^+Cs'$ and so~$w^+$ must lie in~$s'Ct-s'$.
Hence~$w$ must lie in~$s'Ct-t$ (and thus in the interior of~$a_2Ca_1$).
As $a^-_2\in V(A)\subseteq T_3$ and as~$H$ is Hamiltonian connected
by Lemma~\ref{Hlemma1}, there is an $a^-_2$-$w^+$ path~$P$ whose interior consists of
precisely all the vertices in~$H$. But then the $S$-ordered cycle $a_1a_2Cwa^+_1Ca^-_2Pw^+Ca_1$
is Hamiltonian, contradicting the choice of~$C$ (see Figure~\ref{fig:proofcase1}).%
      \COMMENT{We need that $a^+_1$ and $a^-_2$ lie in the interior of $a_1Ca_2$ for
this. But this holds by our choice of~$A_1$ and~$A_2$.}

\medskip

\noindent\textbf{Case~1.2.} \emph{There are no such vertices $a_1\in A_1$ and
$a_2\in A_2$.} 

\smallskip

\noindent
Let~$F_3^-$ denote the set of all predecessors of vertices in~$F_3$.
Recall that~$F$ avoids~$t^+Cs'$. Thus~$F_3^-$ avoids $tCs'-s'$.
Now consider any vertex~$a$ on~$A_2$. Then $N^-_D(a)\subseteq V(C)$
since $a\notin F$ and thus $N^-_D(a)\subseteq V(C)\setminus V(A_1)$
by our assumption. But then using that
$|F_3|\ge (n-k)/2-|H|$ by Lemma~\ref{sizeT5F5}
and arguing similarly as in~(\ref{eq:dega}) one can show that%
     \COMMENT{$d^-_{C-A_1}(a)\ge (n+k)/2-1\ge n-1-|H|-|F_3|$ since
$|F_3|\ge (n-k)/2-|H|$}
$d^-_{C-A_1}(a)\ge n-1-|H|-|F_3|=|C|-1-|F_3|=|C-A_1|-|F_3^-|+|A_1|-1$. Together with the fact that
$F_3^-\cap V(A_1)=\emptyset$ this gives
\begin{equation}\label{eq:daW}
|N^-_{C-A_1}(a)\cap F_3^-|\ge |A_1|-1\ge n/(21k|H|).
\end{equation}
Let~$I_1$ be the subpath of the interval~$I$ preceding the first vertex in~$A_1$.
So $I_1=sCt$. Let $I_2,\dots,I_k$ denote all the other intervals. For each~$i=1,\dots,k$
let~$G_i$ be the auxiliary bipartite graph whose vertex classes are $V(A_2)$ and
$V(I_i)\cap F_3^-$ and in which~$a\in V(A_2)$ is joined to $w\in V(I_i)\cap F_3^-$
if~$wa\in E(D)$. Recall that~$F_3^-$ avoids $tCs'-s'$. Thus
$F_3^-\subseteq V(I_1)\cup\dots\cup V(I_k)$ and so the edges of
$G_1\cup \dots\cup G_k$ correspond to the edges of~$D$ from~$F_3^-$ to~$A_2$.
Together with~(\ref{eq:daW}) this implies that there is some~$i$ such that%
     \COMMENT{Need $n\ge 1260 k^3|H|^2$ for this, ok if $n\ge Ck^9$ where $C$ is suff.~large}
$$
e(G_i)\ge \frac{n|A_2|}{21k^2|H|}\ge \frac{n^2}{420k^3|H|^2}\ge 3n\ge 3|G_i|.
$$
Thus~$G_i$ is not planar and so there are vertices $a_1,a_2\in V(A_2)$
and $w_1,w_2\in V(I_i)\cap F_3^-$ such that the edges $w_1a_1$, $w_2a_2$ `cross' in~$G_i$,
i.e.~such that~$w_1$ lies in the interior of $a_2Cw_2$ and~$a_1$ lies in the
interior of~$w_2Ca_2$.
Recall that $w^+_1,w^+_2\in F_3$ by the definition of~$F_3^-$ and~$a^-_1,a^-_2\in T_3$
as~$A_2$ consisted of inner vertices of~$A$. 
Thus we can apply Corollary~\ref{corH2paths} to obtain disjoint $a^-_j$-$w^+_j$~paths~$P_j$
having all their inner vertices in~$H$ and such that each~$P_j$ contains at least one
inner vertex (where $j=1,2$). Thus $a^-_1P_1w^+_1Cw_2a_2Cw_1a_1Ca^-_2P_2w^+_2Ca^-_1$
is an $S$-ordered cycle%
     \COMMENT{For the cycle to be $S$-ordered we need that~$S$ avoids $w^+_1Cw_2$, but this
holds since $w_1Cw_2\subseteq I_i$ and~$S$ meets~$I_i$ only in its first vertex,
ie $S$ meets $w_1Cw_2$ at most in~$w_1$.}
with at least $|C|+2$ vertices (note that it contains all
the vertices of~$C$), contradicting the choice of~$C$
(see Figure~\ref{fig:proofcase1}).

\medskip

\noindent\textbf{Case~2.} $220k^3\le |H|\le n/2-n/(50k)$.

\smallskip

\noindent
The argument for this case is similar to that in Case~1.
Recall that $|T_3|\ge (n-k)/2-|H| \ge n/(60k)$ by Lemma~\ref{sizeT5F5}
and so one of the~$k$ intervals of~$C$ must contain
at least~$n/(60k^2)$ vertices from~$T_3$. Suppose that this is the case for
the interval~$I$ from~$s$ to~$s'$. Let~$t$ be the first vertex on~$I$ that lies
in~$T_3$. Let~$A$ be the set consisting of the
last~$n/(70k^2)$ vertices from~$T_3$ lying in the interior of~$I$.
For each~$a\in A$ let~$Q_a$ be the set of $220k^3$ vertices of~$C$
preceding~$a$. Note that the definition of~$A$ implies%
    \COMMENT{Need that $n/(60k^2)-n/(70k^2)\ge 220k^3+1$, ok if $n\ge Ck^5$
where $C$ is a suff large constant}
that~$Q_a$ lies in the interior of~$I$ and that~$t$ precedes
the first vertex of~$Q_a$. Together with Lemma~\ref{longloop} this
shows that~$F$ avoids $t^+Cs'$ and thus all of $Q_a\cup \{a,a^+\}$.
In particular, $a\notin F$. Thus $N^-_D(a)\subseteq V(C)$ and so~$a$ satisfies~(\ref{eq:dega}).

\medskip

\noindent\textbf{Case~2.1.} \emph{There is a vertex $a\in A$ for which~$a^+$
receives an edge from some vertex~$q\in Q_a$.} 

\smallskip

\noindent
Inequality~(\ref{eq:dega}) implies that there exists a vertex $w\in N^-_C(a)$
such that the successor~$w^+$ of~$w$ lies in~$F$. Note that~$w$ lies in the interior
of $aCq$ since~$F$ avoids $Q_a\cup \{a,a^+\}$.
As $a\in A\subseteq T_3$ and as~$H$ is Hamiltonian connected by Lemma~\ref{Hlemma1},
there is an $a$-$w^+$ path~$P$ whose interior consists precisely of all the vertices in~$H$.
But then the cycle $qa^+CwaPw^+Cq$ is $S$-ordered and contains $|H|-|Q_a|+1>0$
more vertices than~$C$, a contradiction.

\medskip

\noindent\textbf{Case~2.2.} \emph{There is no such vertex $a\in A$.} 

\smallskip

\noindent
This case is similar to Case~1.2.
Let~$F_3^-$ denote the set of all predecessors of vertices in~$F_3$ again.
Let~$A^+$ denote the set of all successors of vertices in~$A$.
Recall that~$F$ avoids~$t^+Cs'$. Thus~$F_3^-$ avoids $tCs'-s'$ and thus in
particular all the sets~$Q_a$.

Consider any $a\in A$. Then $N^-_D(a^+)\subseteq V(C)$ since $a^+\notin F$
by Corollary~\ref{FTcor}. Thus
$N^-_D(a^+)\subseteq V(C)\setminus Q_a$ by our assumption.
Hence similarly as in Case~1.2 one can show that
$d^-_{C-Q_a}(a^+)\ge |C-Q_a|-|F_3^-|+|Q_a|-1$.
Together with the fact that $F_3^-\cap Q_a=\emptyset$ this gives
\begin{equation}\label{eq:daW2}
|N^-_{C-Q_a}(a^+)\cap F_3^-|\ge |Q_a|-1\ge 210k^3.
\end{equation}
Let~$I_1$ be the subpath of the interval~$I$ preceding the first vertex in~$A^+$.
Let $I_2,\dots,I_k$ denote all the other intervals. For each~$i=1,\dots,k$
let~$G_i$ be the auxiliary bipartite graph whose vertex
classes are~$A^+$ and $V(I_i)\cap F_3^-$ and in which~$a^+\in A^+$ is joined to
$w\in V(I_i)\cap F_3^-$ if~$wa^+$ is an edge of~$D$. Note that
$F_3^-\subseteq V(I_1)\cup\dots\cup V(I_k)$ since~$F_3^-$ avoids $tCs'-s'$.
Thus the edges of $G_1\cup \dots\cup G_k$ correspond to the edges from~$F_3^-$ to~$A^+$.
Together with~(\ref{eq:daW2}) this implies that there is some~$i$ such that
$$
e(G_i)\ge \frac{210k^3|A^+|}{k}=  \frac{210k^3n}{70k^3}= 3n\ge 3|G_i|.
$$
Thus~$G_i$ is not planar and so there are vertices $a^+_1,a^+_2\in V(A^+)$
and $w_1,w_2\in V(I_i)\cap F_3^-$ such that the edges $w_1a^+_1$, $w_2a^+_2$ cross.
As in Case~1.2 we can apply Corollary~\ref{corH2paths} to obtain disjoint $a_j$-$w_j^+$ paths
having all their inner vertices in~$H$ such that each~$P_j$ contains at least one
inner vertex (where $j=1,2$ and~$a_j$ is the predecessor of~$a^+_j$). Thus
$a_1P_1w^+_1Cw_2a^+_2Cw_1a^+_1Ca_2P_2w^+_2Ca_1$
is an $S$-ordered cycle with at least $|C|+2$ vertices (note that it contains
all the vertices of~$C$), contradicting the choice of~$C$.

\medskip

\noindent\textbf{Case~3.} $n/2-n/(50k)\le |H|\le \lceil (n-k)/2\rceil-1$.

\smallskip

\noindent Our first aim is to find vertices $x_1,x_2,y_1,y_2$ on~$C$ with the following
properties:
\begin{enumerate}
\item[(i)] $x_1,x_2,y_1,y_2$ occur on~$C$ in this order and either all of these vertices
are distinct or else $|\{x_1,x_2,y_1,y_2\}|=3$ and $x_1=y_2$.
\item[(ii)] $S$ avoids the interior of $x_1Cx_2$, the interior of~$y_1Cy_2$ as well as~$x_2$
and~$y_1$.
\item[(iii)] There are distinct vertices $h_1,h_2,h'_1,h'_2\in H$ such that
$x_1h_1, x_2h_2, h'_1y_1, h'_2y_2$ are edges.
\item[(iv)] If $x_1\neq y_2$ (and so $x_1,x_2,y_1,y_2$ are distinct) then there are disjoint
$x_i$-$y_i$ paths~$P_i$ of length at least~2 such that all inner vertices of~$P_1$
and~$P_2$ lie in~$H$ and $|P_1\cup P_2|\ge |H|/6$.
\end{enumerate}
To prove the existence of such vertices, suppose first that $|T_8|\ge k+1$
and $|F_8|\ge k+1$. Then we can find two vertices $x_1,x_2\in T_8$ and two
vertices $y_1,y_2\in F_8$ satisfying~(ii). Then these vertices automatically satisfy~(iii). 
Lemma~\ref{longloop} implies that they also satisfy~(i).
Finally, if they are all distinct then Corollary~\ref{corH2paths} shows that they
also satisfy~(iv).

So suppose next that for example $|T_8|\le  k$ but $|F_8|\ge k+1$. Pick $y_1,y_2\in F_8$
as before. To find~$x_1$ and~$x_2$, first note that each vertex $h\in H$ satisfies
$$
d^-_C(h)\ge \delta^-(D)-(|H|-1)\ge \lceil (n+k)/2\rceil -1-  \lceil (n-k)/2\rceil+2
=k+1
$$
and so~$h$ receives at least one edge from some vertex in~$T\setminus T_8$.
As each vertex in~$T\setminus T_8$ sends an edge to at most~7 vertices in~$H$, this means that
there are at least $|H|/7$ independent edges from~$C$ to~$H$.
Thus the interior of some interval of~$C$ contains the endvertices of~16 of these independent edges
which avoid~$y_1$ and~$y_2$. Let~$X_1$ be the set of the first~8 endvertices of these edges
on this interval and let~$X_2$ be the set of the next~8 endvertices. Then Lemma~\ref{H2paths} implies that
there are vertices $x_1\in X_1$ and $x_2\in X_2$ which together with~$y_1$ and~$y_2$
satisfy~(iv). By construction, $x_1,x_2,y_1,y_2$ are all distinct and satisfy~(ii) and~(iii).
Again, Lemma~\ref{longloop} implies that they also satisfy~(i).
The cases when $|T_8|\ge k+1$ but $|F_8|\le k$ and when $|T_8|, |F_8|\le k$ are similar.
So we have shown that there are vertices $x_1,x_2,y_1,y_2$ satisfying (i)--(iv).

In what follows, we will frequently use the fact that any vertex~$x\in V(C)\setminus F_2$ receives
an edge from all but at most
$$|C|-(\delta^-(D)-1)\le n/2+n/(50k)-(n+k)/2+2\le n/(45k)$$
vertices
of~$C$. Similarly, any vertex~$x\in V(C)\setminus T_2$ sends an edge to all but at most~$n/(45k)$
vertices of~$C$.

\medskip

\noindent\textbf{Case~3.1.} $|x_1Cx_2|\ge n/(15k)$

\smallskip

\begin{figure}[htb!]
\begin{center}\footnotesize
\psfrag{1}[][]{$x_1$}
\psfrag{2}[][]{$a_1$}
\psfrag{3}[][]{$A_1$}
\psfrag{4}[][]{$a_2$}
\psfrag{5}[][]{$A_2$}
\psfrag{6}[][]{$x_2$}
\psfrag{7}[][]{$x^+_2$}
\psfrag{8}[][]{$y_1$}
\psfrag{9}[][]{$y^-_2$}
\psfrag{10}[][]{$y_2$}
\psfrag{11}[][]{$P$}
\psfrag{12}[][]{$C$}
\psfrag{13}[][]{$x^+_1$}
\psfrag{14}[][]{$y^-_1$}
\psfrag{15}[][]{$h_2$}
\psfrag{16}[][]{$h'_2$}
\psfrag{17}[][]{$h'_1$}
\psfrag{18}[][]{$h_1$}
\includegraphics[scale=0.55]{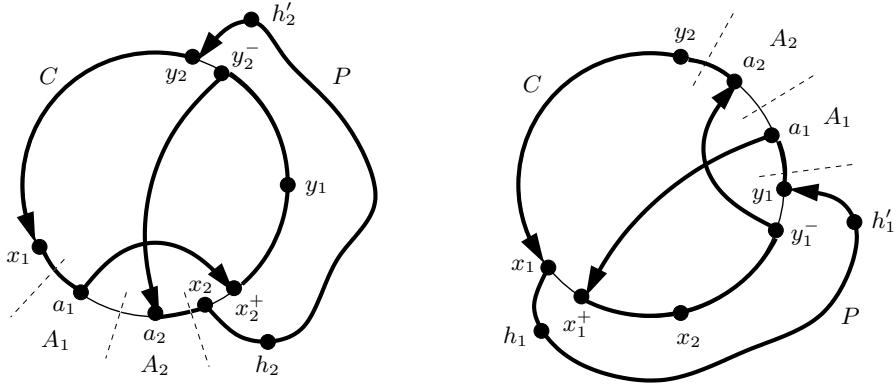}  
\caption{Extending~$C$ into a longer $S$-ordered cycle in Case~3.1 (left) and
Case~3.2 (right).}
\label{fig:proofcase3}
\end{center}
\end{figure}

\noindent
Let~$A_2$ be the set of~$n/(40k)$ vertices which immediately precede~$x_2$ and let~$A_1$
be the set of~$n/(40k)$ vertices which immediately precede~$A_2$. Corollary~\ref{FTcor} implies that
the successor~$x^+_2$ of~$x_2$ on~$C$ does not lie in~$F$. Thus~$x^+_2$ receives an edge
from some vertex $a_1\in A_1$ since it receives an edge from all but at most~$n/(45k)$ vertices of~$C$.
Similarly, the predecessor~$y^-_2$ of~$y_2$ does not lie in~$T$ and thus sends an edge to
some vertex $a_2\in A_2$. Lemma~\ref{Hlemma1} now implies that~$H$ contains a Hamilton path~$P$
from~$h_2$ to~$h'_2$. But then the cycle $a_1x^+_2Cy^-_2a_2Cx_2h_2Ph'_2y_2Ca_1$ is $S$-ordered%
     \COMMENT{We need that $x_2\notin S$ for this.}
and contains all vertices of~$C$ except those in the interior of~$a_1Ca_2$
(see Figure~\ref{fig:proofcase3}).
But as $|H|> n/4>|a_1Ca_2|$ this means that this new cycle is longer than~$C$, a contradiction.

\medskip 

\noindent\textbf{Case~3.2.} $|y_1Cy_2|\ge n/(15k)$

\smallskip

\noindent
The proof of this case is similar to that of Case~3.1.
Let~$A_1$ be the set of~$n/(40k)$ vertices which immediately succeed~$y_1$ and let~$A_2$
be the set of~$n/(40k)$ vertices which immediately succeed~$A_1$. Then the predecessor~$y^-_1$
of~$y_1$ sends an edge to some vertex $a_2\in A_2$ and the successor~$x^+_1$ of~$x_1$
receives an edge from some vertex $a_1\in A_1$. Then the $S$-ordered cycle
$y^-_1a_2Cx_1h_1Ph'_1y_1Ca_1x^+_1Cy^-_1$
is longer%
     \COMMENT{We need that $y_1\notin S$ for the cycle to be $S$-ordered.}
than~$C$, where~$P$ is a Hamilton path in~$H$ from~$h_1$ to~$h'_1$
(see Figure~\ref{fig:proofcase3}).

\medskip

\noindent\textbf{Case~3.3.} $|y_2Cx_1|\ge n/5$

\smallskip

\begin{figure}[htb!]
\begin{center}\footnotesize
\psfrag{1}[][]{$z_1$}
\psfrag{2}[][]{$z_2$}
\psfrag{3}[][]{$z_3$}
\psfrag{4}[][]{$z_4$}
\psfrag{5}[][]{$z_5$}
\psfrag{6}[][]{$z_6$}
\psfrag{7}[][]{$x_1$}
\psfrag{8}[][]{$x^+_1$}
\psfrag{9}[][]{$x_2$}
\psfrag{10}[][]{$x^+_2$}
\psfrag{11}[][]{$y^-_1$}
\psfrag{12}[][]{$y_1$}
\psfrag{13}[][]{$y^-_2$}
\psfrag{14}[][]{$y_2$}
\psfrag{15}[][]{$P_1$}
\psfrag{16}[][]{$P_2$}
\psfrag{17}[][]{$Z_1$}
\psfrag{18}[][]{$Z_2$}
\psfrag{19}[][]{$Z_3$}
\psfrag{20}[][]{$Z_4$}
\psfrag{21}[][]{$Z_5$}
\psfrag{22}[][]{$Z_6$}
\includegraphics[scale=0.55]{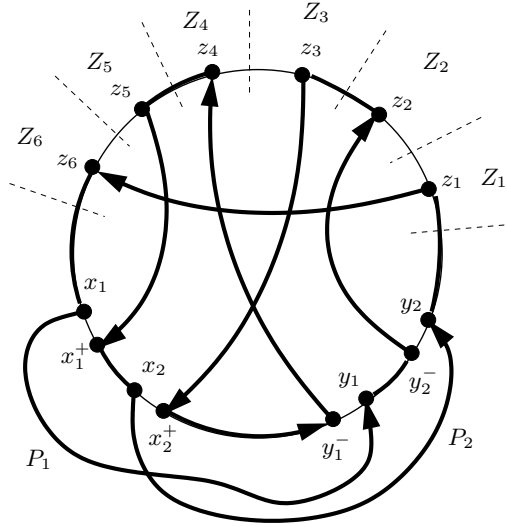}  
\caption{Extending~$C$ into a longer $S$-ordered cycle in Case~3.3.}
\label{fig:proofcase33}
\end{center}
\end{figure}

\noindent
Let~$Z$ be a segment of the interior of $y_2Cx_1$ such that $|Z|\ge n/(6k)$ and such that~$Z$
avoids~$S$. Let~$Z_1$ be the set consisting of the first $n/(40k)$ vertices on~$Z$.
Let~$Z_2$ be the set consisting of the next~$n/(40k)$ vertices and define $Z_3,\dots,Z_6$
similarly. As by Corollary~\ref{FTcor} the predecessor~$y^-_1$ of~$y_1$ does not lie in~$T$
it must send an edge to some vertex $z_4\in Z_4$. Similarly the predecessor~$y^-_2$ of~$y_2$
sends an edge to some vertex $z_2\in Z_2$, the successor~$x^+_1$ of~$x_1$ receives an edge
from some vertex $z_5\in Z_5$ and the successor~$x^+_2$ of~$x_2$ receives an edge
from some vertex $z_3\in Z_3$. Now Lemma~\ref{longloop} implies that either $Z_1\cap T_2=\emptyset$
or $Z_6\cap F_2=\emptyset$ or both. If $Z_1\cap T_2=\emptyset$ then every vertex in~$Z_1$ sends
an edge to~$Z_6$ (since every vertex outside~$T_2$ sends an edge to all but at most~$n/(45k)$
vertices on~$C$). Similarly, if $Z_6\cap F_2=\emptyset$ then every vertex in~$Z_6$ receives
an edge from some vertex in~$Z_1$. So in both cases we can find a $Z_1$-$Z_6$ edge $z_1z_6$.
But then the cycle $x_1P_1y_1Cy^-_2z_2Cz_3x^+_2Cy^-_1z_4Cz_5x^+_1Cx_2P_2y_2Cz_1z_6Cx_1$
is $S$-ordered%
    \COMMENT{we need that $x_2,y_1\notin S$ for this}
and contains at least $|P_1\cup P_2|-4-(|Z|-6)\ge |H|/6-n/(6k)>0$
vertices more than~$C$, a contradiction (see Figure~\ref{fig:proofcase33}).

\medskip

\noindent\textbf{Case~3.4.} \emph{None of Cases~3.1--3.3 holds.}

\smallskip

\begin{figure}[htb!]
\begin{center}\footnotesize
\psfrag{1}[][]{$z_1$}
\psfrag{2}[][]{$z_2$}
\psfrag{3}[][]{$z_3$}
\psfrag{4}[][]{$z_4$}
\psfrag{5}[][]{$z_5$}
\psfrag{6}[][]{$z_6$}
\psfrag{7}[][]{$y^-_1$}
\psfrag{8}[][]{$y_1$}
\psfrag{9}[][]{$y^-_2$}
\psfrag{10}[][]{$y_2$}
\psfrag{11}[][]{$x_1$}
\psfrag{12}[][]{$x^+_1$}
\psfrag{13}[][]{$x_2$}
\psfrag{14}[][]{$x^+_2$}
\psfrag{15}[][]{$P_1$}
\psfrag{16}[][]{$P_2$}
\psfrag{17}[][]{$Z_1$}
\psfrag{18}[][]{$Z_2$}
\psfrag{19}[][]{$Z_3$}
\psfrag{20}[][]{$Z_4$}
\psfrag{21}[][]{$Z_5$}
\psfrag{22}[][]{$Z_6$}
\includegraphics[scale=0.55]{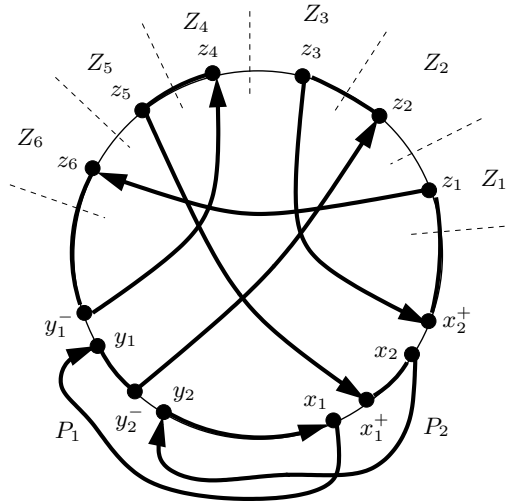}  
\caption{Extending~$C$ into a longer $S$-ordered cycle in Case~3.4.}
\label{fig:proofcase34}
\end{center}
\end{figure}

\noindent
In this case we must have%
     \COMMENT{Indeed, $2n/5+2n/(15k)\le 2n/5+2n/45=20n/45=4n/9<|C|$}
that $|x_2Cy_1|\ge n/5$ and can argue similarly as in
Case~3.3 (see Figure~\ref{fig:proofcase34}). We omit the details.

\medskip

\noindent\textbf{Case~4.} \emph{None of Cases~1--3 holds.}

\smallskip

\noindent
Together with Lemma~\ref{Hlemma1} this implies that~$n-k$ is even
and $|H|=(n-k)/2$. So $|C|=(n+k)/2$.
First note that any vertex $h\in H$ satisfies
\begin{equation}\label{eq:degh}
d^+_C(h), d^-_C(h)\ge (n+k)/2-1-(|H|-1)=k.
\end{equation}
Moreover, if $h,h'\in H$ are distinct and if $s\in S\cap N^-_C(h)$ then by
Lemma~\ref{longloop} the special vertex~$s'$ succeeding~$s$ on~$C$ (i.e.~the unique
vertex $s'\in S$ for which~$S$ avoids the interior of $sCs'$) cannot lie in $N^+_C(h')$.
Thus $|S\cap N^-_C(h)|+|S\cap N^+_C(h')|\le k$ and so
\begin{equation}\label{eq:deghS}
|N^-_C(h)\setminus S|+|N^+_C(h')\setminus S| \ge |N^-_C(h)|+|N^+_C(h')|- k 
\stackrel{(\ref{eq:degh})}{\ge } k.
\end{equation}

\medskip

\noindent\textbf{Case~4.1.} \emph{There exists some vertex $x\in N^-_C(h)\setminus S$.}

\smallskip

\noindent
First note that by Corollary~\ref{FTcor} the successor~$x^+$ of~$x$ on~$C$ does not lie in~$F$.
Thus $d^-_C(x^+)\ge \delta^0(D)= |C|-1$ and so~$x^+$ receives an edge from the predecessor~$x^-$
of~$x$. Pick any vertex $y\in F\setminus \{x,x^-\}$. (Such a vertex exists since $|F|\ge 3$
by~(\ref{eq:degh}).) Note that $y\neq x^+$ since $x^+\notin F$.
By Corollary~\ref{FTcor} the predecessor~$y^-$ of~$y$ does not send an
edge to~$H$ and so $y^-x$ must be an edge (since $d^+_C(y^-)= |C|-1$).
Now apply Lemma~\ref{Hlemma1} to find an $x$-$y$ path~$P$ of length at least~2 all whose
inner vertices lie in~$H$. Then $x^-x^+Cy^-xPyCx^-$ is an $S$-ordered cycle which is longer
than~$C$, a contradiction.

\medskip

\noindent\textbf{Case~4.2.} \emph{There is no vertex as in Case~4.1.}

\smallskip

\noindent
Together with~(\ref{eq:deghS}) this implies that we can find a vertex
$x\in N^+_C(h')\setminus S$. We then argue similarly as in Case~4.1.
This completes the proof of Theorem~\ref{thm:ordHam}.

\section{Acknowledgement}
We are grateful to Oliver Cooley for a careful reading of the manuscript.

\medskip

{\footnotesize \obeylines \parindent=0pt

Daniela K\"{u}hn, Deryk Osthus \& Andrew Young
School of Mathematics
University of Birmingham
Edgbaston
Birmingham
B15 2TT
UK
}

{\footnotesize \parindent=0pt

\it{E-mail addresses}:
\tt{\{kuehn,osthus,younga\}@maths.bham.ac.uk}}

\end{document}